\documentclass[twoside,11pt,leqno]{article}

\newcommand{\be}{\begin{equation}}
\newcommand{\ee}{\end{equation}}
\newcommand{\bea}{\begin{eqnarray}}
\newcommand{\eea}{\end{eqnarray}}
\newcommand{\barray}{\begin{array}}
\newcommand{\earray}{\end{array}}
\newcommand{\pa}{\partial}
\newcommand{\nn}{\nonumber}
\newcommand{\bitem}{\begin{itemize}}
\newcommand{\eitem}{\end{itemize}}
\newtheorem{teo}{Theorem}[section]
\newcommand{\bt}{\begin{teo}}
\newcommand{\et}{\end{teo}}
\newtheorem{Def}{Definition}[section]
\newcommand{\bd}{\begin{Def}}
\newcommand{\ed}{\end{Def}}
\newtheorem{lem}{Lemma}[section]
\newcommand{\bl}{\begin{lem}}
\newcommand{\el}{\end{lem}}
\newtheorem{prop}{Proposition}[section]
\newcommand{\bp}{\begin{prop}}
\newcommand{\ep}{\end{prop}}
\newtheorem{cor}{Corollary}[section]
\newcommand{\bc}{\begin{cor}}
\newcommand{\ec}{\end{cor}}
\newtheorem{ex}{Example}[section]
\newcommand{\bex}{\begin{ex}}
\newcommand{\eex}{\end{ex}}
\newtheorem{rem}{Remark}[section]
\newcommand{\br}{\begin{rem}}
\newcommand{\er}{\end{rem}}

\catcode `\@=11
\@addtoreset{equation}{section}

\begin{document}

\

\bigskip

\begin{center}
{\Large \textbf{COMPATIBLE  FLAT
METRICS\footnote{This work was
supported by the Alexander von Humboldt Foundation (Germany),
the Russian Foundation for Basic Research
(project no. 99--01--00010), and
INTAS (project no. 99--1782).}}}
\end{center}

\smallskip

\begin{center}
{\large {O. I. MOKHOV}}
\end{center}

\medskip

\begin{abstract}
We solve
the problem of description for
nonsingular pairs of compatible flat metrics
in the general $N$-component case.
The integrable nonlinear partial differential equations describing
all nonsingular pairs of compatible flat metrics (or, in other words,
nonsingular flat pencils of metrics) are found and integrated.
The integrating of these equations is based on reducing to
a special nonlinear differential
reduction of the Lam\'{e} equations and
using the Zakharov method of differential
reductions
in the dressing method (a version of the inverse scattering
method).

\end{abstract}

\bigskip

\section{Introduction. Basic definitions} \label{vved}

We shall use both contravariant metrics $g^{ij} (u)$ with upper
indices,
where $u = (u^1,...,u^N)$ are local coordinates, $1 \leq i, j \leq N$,
and covariant metrics
$g_{ij}(u)$ with lower indices,
$g^{is} (u) g_{sj} (u) = \delta^i_j.$
The indices of the
coefficients of the Levi--Civita connections
 $\Gamma^i_{jk} (u)$
and the
indices of the tensors of Riemannian curvature $R^i_{jkl} (u)$
are raised and lowered by the metrics corresponding to them:
\bea
&&
\Gamma^{ij}_k (u) = g^{is} (u) \Gamma^j_{sk} (u),
 \ \ \ \Gamma^i_{jk} (u) = {1 \over 2} g^{is} (u) \left (
{\pa g_{sk} \over \pa u^j} + {\pa g_{js} \over \pa u^k} -
{\pa g_{jk} \over \pa u^s} \right ),\nn\\
&&
R^{ij}_{kl} (u) = g^{is} (u) R^j_{skl} (u), \nn\\
&&
R^i_{jkl} (u) = {\pa \Gamma^i_{jl} \over \pa u^k}
- {\pa \Gamma^i_{jk} \over \pa u^l} +
\Gamma^i_{pk} (u) \Gamma^p_{jl} (u)
- \Gamma^i_{pl} (u) \Gamma^p_{jk} (u).\nn
\eea

\bd \label{d}
Two contravariant flat metrics
$g_1^{ij} (u)$ and $g_2^{ij} (u)$ are called compatible
if any linear combination of these metrics
\be
g^{ij} (u) = \lambda_1 g_1^{ij} (u) + \lambda_2 g_2^{ij} (u),
\label{comb0}
\ee
where $\lambda_1$ and $\lambda_2$ are arbitrary constants such that
$\det ( g^{ij} (u) ) \not\equiv 0$,
is also a flat metric and
the coefficients of the corresponding Levi--Civita connections
are related by the same linear formula:
\be
\Gamma^{ij}_k (u) = \lambda_1 \Gamma^{ij}_{1, k} (u) +
\lambda_2 \Gamma^{ij}_{2, k} (u). \label{sv0}
\ee
In this case,
we shall also say that the flat metrics
$g_1^{ij} (u)$ and $g_2^{ij} (u)$ form a flat pencil
(this definition was proposed by Dubrovin in
{\rm \cite{[31]}, \cite{[30]}}).
\ed

\bd  \label{dd}
Two contravariant metrics
$g_1^{ij} (u)$ and $g_2^{ij} (u)$ of constant Riemannian curvature
$K_1$ and $K_2$, respectively, are called compatible
if any linear combination of these metrics
\be
g^{ij} (u) = \lambda_1 g_1^{ij} (u) + \lambda_2 g_2^{ij} (u),
\label{comb00}
\ee
where $\lambda_1$ and $\lambda_2$ are arbitrary constants such that
$\det ( g^{ij} (u) ) \not\equiv 0$,
is a metric of constant Riemannian curvature
$\lambda_1 K_1 + \lambda_2 K_2$  and
the coefficients of the corresponding Levi--Civita connections
are related by the same linear formula:
\be
\Gamma^{ij}_k (u) = \lambda_1 \Gamma^{ij}_{1, k} (u) +
\lambda_2 \Gamma^{ij}_{2, k} (u). \label{sv00}
\ee
In this case, we shall also say that the metrics
$g_1^{ij} (u)$ and $g_2^{ij} (u)$ form a pencil of
metrics of constant Riemannian curvature.
\ed

\pagestyle{myheadings}
\markboth{{O. I. MOKHOV}\hspace{1.72 in} }
{\hspace{1.15 in} {COMPATIBLE FLAT METRICS}}

\bd  \label{d1}
Two Riemannian or pseudo-Riemannian contravariant metrics
$g_1^{ij} (u)$ and $g_2^{ij} (u)$ are called compatible if
for any linear combination of these metrics
\be
g^{ij} (u) = \lambda_1 g_1^{ij} (u) + \lambda_2 g_2^{ij} (u),
\label{comb}
\ee
where $\lambda_1$ and $\lambda_2$ are arbitrary constants
such that $\det ( g^{ij} (u) ) \not\equiv 0$,
the coefficients of the corresponding Levi--Civita connections
and the components of the corresponding tensors of
Riemannian curvature are related by the same linear formula:
\be
\Gamma^{ij}_k (u) = \lambda_1 \Gamma^{ij}_{1, k} (u) +
\lambda_2 \Gamma^{ij}_{2, k} (u), \label{sv}
\ee
\be
R^{ij}_{kl} (u) = \lambda_1 R^{ij}_{1, kl} (u)
+ \lambda_2 R^{ij}_{2, kl} (u).  \label{kr}
\ee
In this case, we shall also say that the metrics
$g_1^{ij} (u)$ and $g_2^{ij} (u)$ form a pencil of metrics.
\ed

\bd \label{d2}
Two Riemannian or pseudo-Riemannian contravariant metrics
$g_1^{ij} (u)$ and $g_2^{ij} (u)$ are called almost compatible
if for any linear combination of these metrics (\ref{comb})
relation (\ref{sv}) is fulfilled.
\ed

\bd
Two Riemannian or pseudo-Riemannian metrics
$g_1^{ij} (u)$ and $g_2^{ij} (u)$ are called a nonsingular pair
of metrics if the eigenvalues of this pair of metrics, that is,
the roots of the equation
\be
\det ( g_1^{ij} (u) -  \lambda g_2^{ij} (u)) =0,
\ee
are distinct.

A pencil of metrics is called nonsingular if it is formed by
a nonsingular pair of metrics.
\ed

These definitions are motivated by the theory of
compatible Poisson brackets of hydrodynamic type.
We give a brief survey of this theory in the next section.
If the metrics $g_1^{ij} (u)$ and $g_2^{ij}(u)$ are
flat, that is,
 $R^i_{1, jkl} (u) = R^i_{2, jkl} (u) = 0,$ then
relation (\ref{kr}) is equivalent to the condition that
an arbitrary linear combination of the flat metrics
$g_1^{ij} (u)$ and $g_2^{ij}(u)$ is also a flat metric.
In this case, Definition \ref{d1} is equivalent to the well-known
definition of a flat pencil of metrics (Definition \ref{d})
or, in other words,
a compatible pair of local nondegenerate
Poisson structures of hydrodynamic type
\cite{[31]} (see also \cite{[30]}, \cite{[33]},
\cite{[34]}, \cite{[35]}--\cite{mokh1}).
If the metrics $g_1^{ij} (u)$ and $g_2^{ij}(u)$ are
metrics of constant Riemannian curvature $K_1$ and $K_2$,
respectively, that is,
$$R^{ij}_{1, kl} (u) = K_1 (\delta^i_l \delta^j_k -
\delta^i_k  \delta^j_l), \ \ \
R^{ij}_{2, kl} (u) = K_2 (\delta^i_l \delta^j_k -
\delta^i_k  \delta^j_l),$$
then
relation (\ref{kr}) gives the condition that an arbitrary
linear combination of the metrics
 $g_1^{ij} (u)$ and $g_2^{ij}(u)$  (\ref{comb})
is a metric of
constant Riemannian curvature
$\lambda_1 K_1 +
\lambda_2 K_2$. In this case, Definition \ref{d1} is
equivalent to our Definition \ref{dd} of a pencil of
metrics of constant Riemannian curvature or, in other words,
a compatible pair of the corresponding nonlocal Poisson structures
of hydrodynamic type, which were introduced and studied by the
present author
and Ferapontov in \cite{[38]}.
Compatible metrics of more general type correspond
to a compatible pair of nonlocal Poisson structures
of hydrodynamic type that were introduced and studied by
Ferapontov in \cite{[40]}. They arise, for example,
if we use a recursion operator generated by a pair
of compatible Poisson structures of hydrodynamic type.
Such recursion operators
determine, as is well-known, infinite sequences of
corresponding (generally speaking, nonlocal) Poisson structures.

As was earlier noted by the present author in
\cite{[36]}--\cite{mokh1}, condition (\ref{kr}) follows
from condition (\ref{sv}) in the case of certain special
reductions connected with the associativity equations (see also
Theorem \ref{mo1} below).
Of course, it is not accidentally.
Under certain very natural and quite general assumptions on
metrics (it is sufficient but not necessary, in particular,
that eigenvalues of the pair of metrics under
consideration are distinct),
compatibility of the metrics follows from their almost
compatibility but, generally speaking, in the general case,
it is not true even for flat metrics
 (we shall present the corresponding
counterexamples below). Correspondingly,
we would like to emphasize that condition (\ref{sv}),
which is considerably more simple than condition
(\ref{kr}), ``almost'' guarantees compatibility
of metrics and deserves a separate study but, in the general
case, it is necessary to require also
the fulfillment of condition (\ref{kr})
for compatibility of the corresponding Poisson structures of
hydrodynamic type. It is also interesting to find out, does
condition (\ref{kr}) guarantee the fulfillment of condition
(\ref{sv}) or not.

This paper is devoted to the problem of description for
all nonsingular pairs of compatible flat metrics and
to integrability of the corresponding nonlinear partial differential
equations by the inverse scattering method.

\section{Compatible local Poisson structures of
\\ hydrodynamic type}

Any local homogeneous first-order
Poisson bracket, that is,
a Poisson bracket of the form
\be
\{ u^i (x), u^j (y) \} =
g^{ij} (u(x))\, \delta_x (x-y) + b^{ij}_k (u(x)) \, u^k_x \,
\delta (x-y), \label{(1.1)}
\ee
where $u^1,...,u^N$ are local coordinates on a certain smooth
$N$-dimensional manifold $M$, is called a
{\it local Poisson structure of hydrodynamic type} or
{\it Dubrovin--Novikov structure} \cite{[1]}.
Here, $u^i(x),\ 1 \leq i \leq N,$ are functions (fields) of
a single independent variable $x$,
the coefficients $g^{ij}(u)$ and $b^{ij}_k (u)$ of bracket (\ref{(1.1)})
are smooth functions of local coordinates.

In other words,
for arbitrary functionals $I[u]$ and $J[u]$
on the space of fields $u^i(x), \ 1 \leq i \leq N,$ a bracket of the form
\be
\{ I,J \} = \int {\delta I \over \delta u^i(x)}
\biggl ( g^{ij}(u(x)) {d \over dx} + b^{ij}_k (u(x))\, u^k_x \biggr )
{\delta J \over \delta u^j(x)} dx
\label{(1.2)}
\ee
is defined and it is required that this bracket is a Poisson bracket,
that is, it is skew-symmetric:
\be
\{ I, J \} = - \{ J, I \}, \label{skew}
\ee
and satisfies the Jacobi identity
\be
\{ \{ I, J \}, K \} + \{ \{ J, K \}, I \} + \{ \{ K, I \}, J \} =0
\label{jacobi}
\ee
for arbitrary functionals $I[u]$, $J[u]$, and $K[u]$.
The skew-symmetry (\ref{skew}) and the Jacobi identity
(\ref{jacobi}) impose very restrictive conditions on the
coefficients
$g^{ij}(u)$ and $b^{ij}_k (u)$ of bracket (\ref{(1.2)}) (these
conditions will be considered below).
For bracket (\ref{(1.2)}), the Leibniz identity
\be
\{ I \, J, K \} = I \, \{ J, K \} + J \, \{ I, K \}
\label{leib}
\ee
is automatically fulfilled according to the
following property of variational derivative of functionals:
\be
{\delta \left ( I J \right ) \over \delta u^i(x) }
= I {\delta J  \over \delta u^i(x) }
+ J {\delta I \over \delta u^i(x) }.
\ee
Recall that variational derivative of an arbitrary functional
$I [u]$ is defined by
\be
\delta I \equiv I [u + \delta u] - I [u] =
\int {\delta I \over \delta u^k (x)}
\delta u^k (x) dx + o (\delta u).
\ee
The definition of a local Poisson structure of hydrodynamic type
does not depend on a choice of local coordinates $u^1,...,u^N$
on the manifold $M$. Actually, the form of brackets (\ref{(1.2)})
is invariant under local changes of coordinates
$u^i = u^i (v^1,...,v^N),\ 1 \leq i \leq N,$ on $M$:
\bea
&&
\int {\delta I \over \delta u^i(x)}
\biggl ( g^{ij}(u(x)) {d \over dx} + b^{ij}_k (u(x))\, u^k_x \biggr )
{\delta J \over \delta u^j(x)} dx =\\
&&
\int {\delta I \over \delta v^i(x)}
\biggl ( \widetilde g^{ij}(v(x))
{d \over dx} + \widetilde b^{ij}_k (v(x))\, v^k_x \biggr )
{\delta J \over \delta v^j(x)} dx, \nn
\eea
since variational derivatives of functionals transform like
covector fields:
\be
{\delta  I  \over \delta v^i(x) } =
{\delta  I \over \delta u^s(x) } {\pa u^s \over \pa v^i}.
\ee
Correspondingly, the coefficients $g^{ij} (u)$ and
$b^{ij}_k (u)$ of bracket (\ref{(1.2)}) transform as follows:
\be
\widetilde g^{sr} (v) = g^{ij} (u(v)){\partial v^s \over \partial u^i}
 {\partial v^r \over \partial u^j},
\label{(1.3)}
\ee

\be
\widetilde b^{sr}_l (v) = b^{ij}_k (u(v)) {\partial v^s \over \partial u^i}
{\partial v^r \over \partial u^j}
 {\partial u^k \over \partial v^l} + g^{ij}(u(v))
{\partial v^s \over \partial u^i}
{\partial^2 v^r \over \partial u^j \partial u^p}
{\partial u^p \over \partial v^l}.
\label{(1.4)}
\ee
In particular, the coefficients $g^{ij}(u)$
define a contravariant tensor field of rank 2
(a contravariant ``metric'') on the manifold $M$.
For the important case of a nondegenerate metric
$g^{ij} (u), \ \det g^{ij} \neq 0,$ (that is,
in the case of a pseudo-Riemannian manifold
$(M, g^{ij})$), the coefficients
$b^{ij}_k (u)$ define the Christoffel symbols
of an affine connection
$\Gamma^i_{jk} (u)$:
\be
b^{ij}_k (u) = - g^{is} (u) \Gamma^j_{sk} (u), \label{kris}
\ee
$$\widetilde \Gamma^i_{jk} (v) = \Gamma^p_{rs} (u(v))
{\partial v^i \over \partial u^p}
{\partial u^r \over \partial v^j}
{\partial u^s \over \partial v^k} +
{\partial^2 u^s \over \partial v^j \partial v^k}
{\partial v^i \over \partial u^s}. $$

The local Poisson structures of hydrodynamic type (\ref{(1.1)})
were introduced and studied by Dubrovin and Novikov
in \cite{[1]}. In this paper, they proposed
a general local Hamiltonian approach (this approach
corresponds to the local structures of form (\ref{(1.1)})) to the so-called
{\it homogeneous systems of hydrodynamic type}, that is,
evolutionary quasilinear systems of first-order partial differential
equations
\be
u^i_t = V^i_j (u)\, u^j_x.
\label{(1.5)}
\ee

This Hamiltonian approach was motivated by the study
of the equations of Euler hydrodynamics and the
Whitham averaging equations, which describe the evolution of
slowly modulated multiphase solutions of partial
differential equations \cite{[3]}.

Local bracket (\ref{(1.2)}) is called
{\it nondegenerate} if
$\det (g^{ij} (u)) \not\equiv 0$.
For the general nondegenerate brackets of form (\ref{(1.2)}),
Dubrovin and Novikov proved the following important theorem.

\bt [Dubrovin and Novikov \cite{[1]}] \label{dn}
If $\det (g^{ij} (u)) \not\equiv 0$, then bracket (\ref{(1.2)})
is a Poisson bracket, that is, it is skew-symmetric and satisfies
the Jacobi identity, if and only if
\bitem
\item [(1)] $g^{ij} (u)$ is an arbitrary flat
pseudo-Riemannian contravariant metric (a metric of
zero Riemannian curvature),

\item [(2)] $b^{ij}_k (u) = - g^{is} (u) \Gamma ^j_{sk} (u),$ where
$\Gamma^j_{sk} (u)$ is the Riemannian connection generated by
the contravariant metric $g^{ij} (u)$
(the Levi--Civita connection).
\eitem
\et

Consequently, for any local nondegenerate Poisson structure of
hydrodynamic type, there always exist local coordinates
$v^1,...,v^N$ (flat coordinates of the metric $g^{ij}(u)$) in which
all the coefficients of the bracket are constant:
\be
\widetilde g^{ij} (v)
= \eta^{ij} = {\rm \ const}, \ \
\widetilde \Gamma^i_{jk} (v) = 0, \ \
\widetilde b^{ij}_k (v) =0,
\ee
that is, the bracket has the constant form
\be
\{ I,J \} = \int {\delta I \over \delta v^i(x)}
 \eta^{ij} {d \over dx}
{\delta J \over \delta v^j(x)} dx,
\label{(1.6)}
\ee
where $(\eta^{ij})$ is a nondegenerate symmetric constant matrix:
\be
\eta^{ij} = \eta^{ji}, \ \ \eta^{ij} = {\rm const},
\  \ \det \, (\eta^{ij}) \neq 0.\nn
\ee

On the other hand, as early as 1978, Magri proposed
a bi-Hamil\-to\-ni\-an approach to the integration of nonlinear
systems \cite{[11]}. This approach demonstrated
that integrability is closely related to the
bi-Ha\-mil\-to\-ni\-an property, that is, to the property
of a system to have two compatible Hamiltonian
representations. As was shown by Magri in \cite{[11]},
compatible Poisson brackets generate integrable hierarchies of
systems of differential equations.
Therefore, the description of compatible
Poisson structures is very urgent
and important problem in the theory of
integrable systems.
In particular, for a system, the bi-Hamiltonian property
generates recurrent relations for the conservation laws of
this system.

Beginning from \cite{[11]}, quite extensive literature
(see, for example,
\cite{[16]}, \cite{[14]}, \cite{[13]}, \cite{[12]},
\cite{[15]},
and the necessary references therein) has been
devoted to the bi-Hamil\-to\-ni\-an approach and to
the construction of compatible Poisson structures
for many specific important equations of
mathematical physics and field theory.
As far as the problem of description
of sufficiently wide classes
of compatible Poisson structures of defined special types is
concerned, apparently the first such statement was considered
in \cite{[17]}, \cite{[18]} (see also \cite{[19]}, \cite{[20]}).
In those papers, the present author posed and completely solved
the problem of description of all compatible local scalar first-order
and third-order Poisson brackets, that is, all Poisson brackets given
by arbitrary scalar first-order and third-order ordinary
differential
operators.
These brackets generalize the well-known compatible pair of
the Gardner--Zakharov--Faddeev bracket
\cite{[21]}, \cite{[22]}
(the first-order bracket) and the Magri bracket \cite{[11]}
(the third-order bracket) for the Korteweg--de Vries equation.

In the case of homogeneous systems of hydrodynamic type, many integrable
systems possess compatible Poisson structures of hydrodynamic type.
The problems of description of these structures for particular systems and
numerous examples were considered in many papers
(see, for example, \cite{[23]}, \cite{[24]}, \cite{[25]},
\cite{[26]}, \cite{[27]}, \cite{[29]}).
In particular, in \cite{[23]}
Nutku studied a special class of compatible
two-component Poisson structures of hydrodynamic type and
the related bi-Hamiltonian hydrodynamic systems.
In \cite{[28]} Ferapontov classified all two-component
homogeneous systems of hydrodynamic type possessing
three compatible nondegenerate local Poisson structures of
hydrodynamic type.

In the general form, the problem of description of
flat pencils of metrics (or, in other words,
compatible nondegenerate local Poisson structures of hydrodynamic type)
was considered by Dubrovin in
\cite{[31]}, \cite{[30]}
in connection with the construction of important examples of
such flat pencils of metrics, generated by natural pairs of
flat metrics on the spaces of orbits of Coxeter groups and on other
Frobenius manifolds and associated with the
corresponding quasi-homogeneous solutions of
the associativity equations.
In the theory of Frobenius manifolds introduced and studied by
Dubrovin
\cite{[31]}, \cite{[30]}
(they correspond to two-dimensional topological field theories),
a key role is played by flat pencils of metrics, possessing
a number of special additional (and very restrictive) properties
(they satisfy the so-called quasi-homogeneity property).
In addition, in \cite{[33]} Dubrovin proved that the theory of
Frobenius manifolds is equivalent to the theory of
quasi-homogeneous compatible nondegenerate local Poisson structures of
hydrodynamic type. The general problem on compatible
nondegenerate local Poisson structures of hydrodynamic type
was also considered by Ferapontov
in \cite{[34]}.

The present author's papers
\cite{[35]}--\cite{mokh1}
are devoted to the general problem of classification of
all compatible local Poisson structures of hydrodynamic type,
the study of the integrable
nonlinear systems that describe such
the compatible Poisson structures and, mainly,
the special reductions connected
with the associativity equations.

\bd [Magri \cite{[11]}] \label{dm}
Two Poisson brackets $\{ \ , \ \}_1$
and $\{ \ , \ \}_2$ are called {\it compatible} if
an arbitrary linear combination of these Poisson brackets
\be
\{ \ , \ \} = \lambda_1 \, \{ \ , \ \}_1 +
\lambda_2 \, \{ \ , \ \}_2, \label{magri}
\ee
where $\lambda_1$ and $\lambda_2$ are arbitrary constants,
is also a Poisson bracket.
In this case, we also say that the brackets $\{ \ , \ \}_1$ and
$\{ \ , \ \}_2$ form a pencil of Poisson brackets.
\ed

Correspondingly, the problem of description
for compatible nondegenerate local Poisson structures
of hydrodynamic type is pure diffe\-ren\-ti\-al-geometric
problem of description for flat pencils of metrics
(see \cite{[31]}, \cite{[30]}).

In \cite{[31]}, \cite{[30]}
Dubrovin presented all the tensor relations
for the general flat pencils of metrics.
First, we introduce the necessary notation.
Let $\nabla_1$ and $\nabla_2$
be the operators of covariant differentiation given by the Levi--Civita
connections $\Gamma^{ij}_{1,k} (u)$
and $\Gamma^{ij}_{2,k} (u)$, generated by the metrics $g^{ij}_1 (u)$
and $g^{ij}_2 (u)$, respectively. The indices
of the covariant differentials are raised and lowered by
the
corresponding metrics: $\nabla^i_1= g^{is}_1 (u) \nabla_{1,s}$,
$\nabla^i_2=g^{is}_2 (u) \nabla_{2,s}$.
Consider the tensor
\be
\Delta ^{ijk} (u) = g^{is}_1 (u) g^{jp}_2 (u)
\left (\Gamma^k_{2, ps} (u)
- \Gamma^k_{1, ps} (u) \right )  \label{(2.3)}
\ee
introduced by Dubrovin in \cite{[31]}, \cite{[30]}.

\bt [Dubrovin \cite{[31]}, \cite{[30]}] \label{dub1}
 If metrics
$g^{ij}_1 (u)$ and $g^{ij}_2 (u)$ form a flat pencil,
then there exists a vector field $f^i (u)$ such that the tensor
$\Delta ^{ijk} (u)$ and the metric $g^{ij}_1 (u)$ have the form
\be
\Delta ^{ijk} (u) = \nabla^i_2 \nabla^j_2 f^k (u),
\label{(2.4)}
\ee
\be
g^{ij}_1 (u) = \nabla^i_2 f^j (u) + \nabla^j_2 f^i (u) +
c g^{ij}_2 (u), \label{(2.5)}
\ee
where $c$ is a certain constant,
and the vector field $f^i (u)$ satisfies the equations
\be
\Delta^{ij}_s (u) \Delta^{sk}_l (u) =
\Delta^{ik}_s (u) \Delta^{sj}_l (u), \label{(2.6)}
\ee
where
\be
\Delta^{ij}_k (u) =g_{2,ks} (u) \Delta ^{sij} (u)
= \nabla_{2,k} \nabla^i_2 f^j (u), \label{(2.7)}
\ee
and
\be
(g^{is}_1 (u) g^{jp}_2 (u) - g^{is}_2 (u) g^{jp}_1 (u))
\nabla_{2,s} \nabla_{2,p} f^k (u) =0. \label{(2.8)}
\ee
Conversely, for the flat metric $g^{ij}_2 (u)$ and the vector field
$f^i (u)$ that is a solution of the system of equations
(\ref{(2.6)}) and (\ref{(2.8)}),
the metrics $g^{ij}_2 (u)$ and (\ref{(2.5)}) form a flat pencil.
\et

The proof of this theorem immediately follows from the relations that
are equivalent to the fact that the metrics $g^{ij}_1 (u)$
and $g^{ij}_2 (u)$ form a flat pencil and are considered in
flat coordinates of the metric $g^{ij}_2 (u)$ \cite{[31]}, \cite{[30]}.

In my paper \cite{[35]}, an explicit and simple {\em criterion}
of compatibility for two local Poisson structures of
hydrodynamic type is formulated, that is, it is shown
what explicit form is sufficient and necessary for the local
Poisson structures of hydrodynamic type to be
compatible.

For the moment, in the general case, we are able to formulate such an
explicit criterion
only namely in terms of Poisson structures but not in terms of
metrics as in Theorem \ref{dub1}. But for nonsingular pairs of
the Poisson structures of hydrodynamic type
(that is, for nonsingular pairs of the corresponding
metrics), in this paper, we shall get
an explicit general criterion of compatibility namely
in terms of the corresponding metrics.

\bl [\cite{[35]}]  \label{lem1}
{\bf (An explicit criterion of
compatibility for local Poisson structures of
hydrodynamic type)}
Any local Poisson structure of hydrodynamic type
$\{ I, J \}_2$  is compatible with the constant nondegenerate
Poisson bracket (\ref{(1.6)}) if and only if it has the form
\bea
&&
\{ I, J \}_2 =\label{(2.9)}\\
&&
\int {\delta I \over \delta v^i (x)} \biggl (
\biggl ( \eta^{is} {\partial h^j \over \partial v^s}
+ \eta^{js} {\partial h^i \over \partial v^s} \biggr )
{d \over dx} + \eta^{is}
{\partial^2 h^j \over \partial v^s \partial v^k}
v^k_x \biggr ) {\delta J \over \delta v^j (x)} dx, \nn
\eea
where $h^i (v), \ 1 \leq i \leq N,$ are smooth
functions defined on a certain neighbourhood.
\el

We do not require in Lemma \ref{lem1} that
the Poisson structure of hydrodynamic type
$\{ I, J \}_2$ is nondegenerate.
Besides, it is important to note that this statement is local.

In 1995, in the paper \cite{[34]}, Ferapontov
proposed an approach to the problem on flat pencils of metrics,
which is motivated by the theory of recursion operators,
and formulated (without any proof) the following theorem as a criterion of
compatibility for nondegenerate local Poisson structures
of hydrodynamic type:

\bt [\cite{[34]}] \label{tfer}
Two local nondegenerate Poisson structures of hydrodynamic type
given by flat metrics $g_1^{ij}(u)$ and $g_2^{ij}(u)$
are compatible if and only if the Nijenhuis tensor of the affinor
$v^i_j (u) = g_1^{is} (u) g_{2, sj} (u)$ vanishes,
that is,
\be
\ \ \ N^k_{ij} (u) = v^s_i (u) {\pa v^k_j \over \pa u^s}
- v^s_j (u) {\pa v^k_i \over \pa u^s}  +
v^k_s (u) {\pa v^s_i \over \pa u^j} -
v^k_s (u) {\pa v^s_j \over \pa u^i} =0.
\ee
\et

Besides, in the remark in \cite{[34]}, it is noted that
if the spectrum of $v^i_j (u)$ is simple,
then the vanishing of the Nijenhuis tensor implies
the existence of coordinates $R^1,...,R^N$ for which
all the objects $v^i_j (u)$, $g_1^{ij} (u)$, $g_2^{ij} (u)$
become diagonal.
Moreover, in these coordinates the
$i$th eigenvalue of $v^i_j (u)$ depends only on the
coordinate $R^i$.
In the case when all the eigenvalues are nonconstant,
they can be introduced as new coordinates.
In these new coordinates
$\tilde v^i_j (R) = {\rm diag \ } (R^1,...,R^N)$,
$\tilde g_2^{ij} (R) = {\rm diag \ } (g^1 (R),...,g^N (R))$,
$\tilde g_1^{ij} (R) = {\rm diag \ } (R^1 g^1 (R),..., R^N g^N (R))$.

In this paper, we shall prove that, unfortunately,
in the general case, Theorem \ref{tfer} is not true
and, correspondingly, it is not a criterion of
compatibility of flat metrics.
Generally speaking, compatibility of flat metrics does not follow
from the vanishing of the corresponding Nijenhuis tensor.
The corresponding counterexamples will be presented below
in Section \ref{counter}.
We also prove that,
in the general case, Theorem \ref{tfer} is actually a criterion of
almost compatibility of flat metrics that
does not guarantee compatibility of the corresponding
nondegenerate local Poisson structures of hydrodynamic type.
But if the spectrum of $v^i_j (u)$ is simple, that is, all the
eigenvalues are distinct, then we prove that
Theorem \ref{tfer} is not only true but also can be
essentially  generalized to the case of arbitrary
compatible Riemannian or pseudo-Riemannian metrics,
in particular, the especially important cases in
the theory of systems of hydrodynamic type, namely,
the cases of metrics of constant Riemannian curvature
or the metrics generating the general nonlocal Poisson
structures of hydrodynamic type.

Namely, we shall prove
the following theorems for any pseudo-Rie\-man\-ni\-an metrics
(not only for flat metrics as in Theorem \ref{tfer}).
\bt \label{tmo1}
\bitem
\item[1)] If for any linear combination (\ref{comb})
of two metrics $g_1^{ij} (u)$ and $g_2^{ij} (u)$ condition (\ref{sv})
is fulfilled, then the Nijenhuis tensor of the affinor
$$v^i_j (u) = g_1^{is} (u) g_{2, sj} (u)$$
vanishes.
Thus, for any two compatible or almost compatible metrics,
the corresponding Nijenhuis tensor always vanishes.
\item[2)] If a pair of metrics
$g_1^{ij} (u)$ and $g_2^{ij} (u)$ is nonsingular, that is,
the roots of the equation
\be
\det ( g_1^{ij} (u) -  \lambda g_2^{ij} (u)) =0
\ee
are distinct, then it follows from the vanishing
of the Nijenhuis tensor of the affinor
$v^i_j (u) = g_1^{is} (u) g_{2, sj} (u)$ that the metrics
$g_1^{ij} (u)$ and $g_2^{ij} (u)$ are compatible.
Thus, a nonsingular pair of metrics is compatible
if and only if the metrics are almost compatible.
\eitem
\et
\bt    \label{teomxx}
Any nonsingular pair of metrics is compatible
if and only if there exist local coordinates $u = (u^1,...,u^N)$
such that
$g^{ij}_2 (u) = g^i (u) \delta^{ij}$ and
$g^{ij}_1 (u) = f^i (u^i) g^i (u) \delta^{ij},$
where $f^i (u^i),$ $i=1,...,N,$ are arbitrary (generally speaking,
complex)
functions
of single variables (of course,
the functions $f^i (u^i)$
are not equal identically to zero and,
 for nonsingular pairs of metrics,
all these functions must be distinct and they can not be equal
to one another if they are constants but,
nevertheless, in this special case, the metrics will be also compatible).
\et

Sections 3 and 4 are devoted to the proof
of Theorems \ref{tmo1} and \ref{teomxx}.

\section{Almost compatible metrics and \\
the Nijenhuis tensor} \label{sect1}

Let us consider two arbitrary
contravariant Riemannian or pseudo-Riemannian
metrics
$g_1^{ij} (u)$  and  $g_2^{ij} (u)$, and also
the corresponding coefficients of the Levi--Civita connections
$\Gamma^{ij}_{1, k} (u)$ and $\Gamma^{ij}_{2, k} (u)$.

We introduce the tensor

\bea
&&
M^{ijk} (u) = \label{tensorm}\\
&&
g_1^{is} (u) \Gamma_{2, s}^{jk} (u) -
g_2^{js} (u) \Gamma_{1, s}^{ik} (u) -
g_1^{js} (u) \Gamma_{2, s}^{ik} (u) +
g_2^{is} (u) \Gamma_{1, s}^{jk} (u).\nn
\eea

It follows from the following representation
that $M^{ijk} (u)$ is actually a tensor:
\bea
&&
M^{ijk} (u) = g_1^{is} (u) g_2^{jp} (u) ( \Gamma_{2, ps}^k (u) -
 \Gamma_{1, ps}^k (u) ) -\\
&&
g_1^{js} (u) g_2^{ip} (u) (\Gamma_{2, ps}^k (u) -
 \Gamma_{1, ps}^k (u)).\nn
\eea

\bl \label{l1}
The tensor $M^{ijk} (u)$ vanishes if and only if
the metrics $g_1^{ij} (u)$ and $g_2^{ij} (u)$ are almost
compatible.
\el

{\it Proof}. Recall that functions $\Gamma^{ij}_k (u)$ define
the Christoffel symbols of the Levi--Civita connection
for a contravariant metric $g^{ij} (u)$
if and only if the following relations are fulfilled:
\be
{\pa g^{ij} \over \pa u^k} + \Gamma^{ij}_k (u)
+ \Gamma^{ji}_k (u) =0,\label{ch1}
\ee
that is, the connection is compatible with the metric, and
\be
g^{is} (u) \Gamma^{jk}_s (u) = g^{js} (u) \Gamma^{ik}_s (u), \label{ch2}
\ee
that is, the connection is symmetric.

If $g^{ij} (u)$ and $\Gamma^{ij}_k (u)$ are defined by
formulas (\ref{comb}) and (\ref{sv}), respectively,
 then linear relation (\ref{ch1})
is automatically fulfilled and relation (\ref{ch2}) is exactly
equivalent to the relation $M^{ijk} (u) =0$. q.e.d.

Let us introduce the affinor
\be
v^i_j (u) = g_1^{is} (u) g_{2, sj} (u) \label{aff}
\ee
and consider the Nijenhuis tensor of this affinor
\be
N^k_{ij} (u) = v^s_i (u) {\pa v^k_j \over \pa u^s}
- v^s_j (u) {\pa v^k_i \over \pa u^s}  +
v^k_s (u) {\pa v^s_i \over \pa u^j} -
v^k_s (u) {\pa v^s_j \over \pa u^i},  \label{nij}
\ee
following \cite{[34]}, where were similarly considered
the affinor $v^i_j (u)$
and its Nijenhuis tensor for two flat metrics.

\bt  \label{pochti}
Any two metrics
$g_1^{ij} (u)$ and $g_2^{ij} (u)$ are almost
compatible if and only if the corresponding Nijenhuis tensor
$N^k_{ij} (u)$ (\ref{nij}) vanishes.
\et

\bl \label{l2}
The following identities are always fulfilled:
\bea
&&
g_{1, sp} (u) N^p_{rq} (u) g_2^{ri} (u) g_2^{qj} (u) g_2^{sk} (u)
=\label{mn1}\\
&&
M^{kji} (u) + M^{ikj} (u) + M^{ijk} (u),\nn
\eea
\bea
&&
2( M^{ikj} (u) + M^{ijk} (u))=
g_{1, sp} (u) N^p_{rq} (u) g_2^{ri} (u) g_2^{qj} (u) g_2^{sk} (u) +
\label{mn2}\\
&&
g_{1, sp} (u) N^p_{rq} (u) g_2^{ri} (u) g_2^{qk} (u) g_2^{sj} (u),\nn
\eea
\bea
&&
2 M^{kji} (u) =
g_{1, sp} (u) N^p_{rq} (u) g_2^{ri} (u) g_2^{qj} (u) g_2^{sk} (u) -
\label{mn3}\\
&&
g_{1, sp} (u) N^p_{rq} (u) g_2^{ri} (u) g_2^{qk} (u) g_2^{sj} (u).\nn
\eea
\el

{\it Proof.}  In the following
calculations, using many times
both relations (\ref{ch1}) and
(\ref{ch2}) for both the metrics $g^{ij}_1 (u)$
and $g^{ij}_2 (u)$, we have
\bea
&&
N^k_{ij} (u) = v^s_i  {\pa v^k_j \over \pa u^s}
- v^s_j  {\pa v^k_i \over \pa u^s}  +
v^k_s {\pa v^s_i \over \pa u^j} -
v^k_s  {\pa v^s_j \over \pa u^i} =
\nn\\
&&
g^{sp}_1 g_{2, pi} {\pa \over \pa u^s} \left (
g^{kl}_1 g_{2, lj} \right ) -
g^{sp}_1 g_{2, pj} {\pa \over \pa u^s} \left (
g^{kl}_1 g_{2, li} \right ) +\nn\\
&&
g^{kp}_1 g_{2, ps} {\pa \over \pa u^j} \left (
g^{sl}_1 g_{2, li} \right ) -
g^{kp}_1 g_{2, ps} {\pa \over \pa u^i} \left (
g^{sl}_1 g_{2, lj} \right ) =\nn\\
&&
- g^{sp}_1 g_{2, pi} g_{2, lj} \left (
\Gamma^{kl}_{1, s} + \Gamma^{lk}_{1, s} \right ) +
g^{sp}_1 g_{2, pi} g^{kl}_1 g_{2, lr} g_{2, tj} \left (
\Gamma^{rt}_{2, s} + \Gamma^{tr}_{2, s} \right ) +\nn\\
&&
g^{sp}_1 g_{2, pj} g_{2, li} \left (
\Gamma^{kl}_{1, s} + \Gamma^{lk}_{1, s} \right ) -
g^{sp}_1 g_{2, pj} g^{kl}_1 g_{2, lr} g_{2, ti} \left (
\Gamma^{rt}_{2, s} + \Gamma^{tr}_{2, s} \right ) -\nn\\
&&
g^{kp}_1 g_{2, ps} g_{2, li} \left (
\Gamma^{sl}_{1, j} + \Gamma^{ls}_{1, j} \right ) +
g^{kp}_1 g_{2, ps} g^{sl}_1 g_{2, lr} g_{2, ti} \left (
\Gamma^{rt}_{2, j} + \Gamma^{tr}_{2, j} \right ) +\nn\\
&&
g^{kp}_1 g_{2, ps} g_{2, lj} \left (
\Gamma^{sl}_{1, i} + \Gamma^{ls}_{1, i} \right ) -
g^{kp}_1 g_{2, ps} g^{sl}_1 g_{2, lr} g_{2, tj} \left (
\Gamma^{rt}_{2, i} + \Gamma^{tr}_{2, i} \right ),\nn
\eea

\bea
&&
N^k_{ij} g^{in}_2 g^{jm}_2=\nn\\
&&
- g^{sn}_1 \left (
\Gamma^{km}_{1, s} + \Gamma^{mk}_{1, s} \right ) +
g^{sn}_1 g^{kl}_1 g_{2, lr} \left (
\Gamma^{rm}_{2, s} + \Gamma^{mr}_{2, s} \right ) +\nn\\
&&
g^{sm}_1  \left (
\Gamma^{kn}_{1, s} + \Gamma^{nk}_{1, s} \right ) -
g^{sm}_1 g^{kl}_1 g_{2, lr} \left (
\Gamma^{rn}_{2, s} + \Gamma^{nr}_{2, s} \right ) -\nn\\
&&
g^{kp}_1 g_{2, ps} g^{jm}_2 \left (
\Gamma^{sn}_{1, j} + \Gamma^{ns}_{1, j} \right ) +
g^{kp}_1 g_{2, ps} g^{sl}_1 g_{2, lr} g^{jm}_2 \left (
\Gamma^{rn}_{2, j} + \Gamma^{nr}_{2, j} \right ) +\nn\\
&&
g^{kp}_1 g_{2, ps} g^{in}_2 \left (
\Gamma^{sm}_{1, i} + \Gamma^{ms}_{1, i} \right ) -
g^{kp}_1 g_{2, ps} g^{sl}_1 g_{2, lr} g^{in}_2 \left (
\Gamma^{rm}_{2, i} + \Gamma^{mr}_{2, i} \right )=\nn\\
&&
- g^{sn}_1
\Gamma^{km}_{1, s} +
g^{sn}_1 g^{kl}_1 g_{2, lr} \left (
\Gamma^{rm}_{2, s} + \Gamma^{mr}_{2, s} \right ) +\nn\\
&&
g^{sm}_1
\Gamma^{kn}_{1, s} -
g^{sm}_1 g^{kl}_1 g_{2, lr} \left (
\Gamma^{rn}_{2, s} + \Gamma^{nr}_{2, s} \right ) -\nn\\
&&
g^{kp}_1 g_{2, ps} g^{jm}_2 \left (
\Gamma^{sn}_{1, j} + \Gamma^{ns}_{1, j} \right ) +
g^{kp}_1 g_{2, ps} g^{sj}_1
\Gamma^{mn}_{2, j} +\nn\\
&&
g^{kp}_1 g_{2, ps} g^{in}_2 \left (
\Gamma^{sm}_{1, i} + \Gamma^{ms}_{1, i} \right ) -
g^{kp}_1 g_{2, ps} g^{si}_1
\Gamma^{nm}_{2, i},\nn
\eea

\bea
&&
g_{1, qk} N^k_{ij} g^{in}_2 g^{jm}_2 =
- \Gamma^{nm}_{1, q} +
g^{sn}_1 g_{2, qr} \left (
\Gamma^{rm}_{2, s} + \Gamma^{mr}_{2, s} \right ) +\nn\\
&&
\Gamma^{mn}_{1, q} -
g^{sm}_1 g_{2, qr} \left (
\Gamma^{rn}_{2, s} + \Gamma^{nr}_{2, s} \right ) -
g_{2, qs} g^{jm}_2 \left (
\Gamma^{sn}_{1, j} + \Gamma^{ns}_{1, j} \right ) +\nn\\
&&
g_{2, qs} g^{sj}_1
\Gamma^{mn}_{2, j} +
g_{2, qs} g^{in}_2 \left (
\Gamma^{sm}_{1, i} + \Gamma^{ms}_{1, i} \right ) -
g_{2, qs} g^{si}_1
\Gamma^{nm}_{2, i},\nn
\eea
and, finally,
\bea
&&
g_{1, qk} N^k_{ij} g^{in}_2 g^{jm}_2 g^{tq}_2 =
- g^{tq}_2 \Gamma^{nm}_{1, q} +
g^{sn}_1 \left (
\Gamma^{tm}_{2, s} + \Gamma^{mt}_{2, s} \right ) +\nn\\
&&
g^{tq}_2 \Gamma^{mn}_{1, q} -
g^{sm}_1 \left (
\Gamma^{tn}_{2, s} + \Gamma^{nt}_{2, s} \right ) -
 g^{jm}_2 \left (
\Gamma^{tn}_{1, j} + \Gamma^{nt}_{1, j} \right ) +\nn\\
&&
g^{tj}_1
\Gamma^{mn}_{2, j} +
g^{in}_2 \left (
\Gamma^{tm}_{1, i} + \Gamma^{mt}_{1, i} \right ) -
g^{ti}_1
\Gamma^{nm}_{2, i}=\nn\\
&&
M^{tmn} + M^{ntm} + M^{nmt}.\nn
\eea

Note that the tensor $M^{ijk} (u)$ (\ref{tensorm})
is skew-symmetric
with respect to the indices $i$ and $j$.
Permuting the indices  $k$ and $j$ in formula (\ref{mn1}) and
adding the corresponding relation to (\ref{mn1}),
we obtain (\ref{mn2}).
Formula (\ref{mn3}) follows from
(\ref{mn1}) and (\ref{mn2}) straightforward. q.e.d.

\bc
The tensor $M^{ijk} (u)$ vanishes if and only if
the Nijenhuis tensor (\ref{nij}) vanishes.
\ec

In the papers \cite{[35]}--\cite{mokh1}, the present author
studied special
reductions in the general problem on compatible flat metrics,
namely, the reductions
connected with the associativity equations, that is, the following
general ansatz in formula (\ref{(2.9)}):
$$
h^i (v) = \eta^{is} {\pa \Phi \over \pa v^s},
$$
where
$\Phi (v^1,...,v^N)$ is a function of $N$ variables.

Correspondingly,
in this case the metrics have the form:

\be
g_1^{ij} (v) = \eta^{ij}, \ \ \ g_2^{ij} (v) = \eta^{is}
\eta^{jp} {\pa^2 \Phi \over \pa v^s \pa v^p}.  \label{m1}
\ee

\bt [\cite{[36]}, \cite{mokh2}, \cite{mokh1}] \label{mo1}
If metrics (\ref{m1}) are almost compatible, then
they are compatible.
Moreover, in this case, the metric $g_2^{ij} (v)$ is
necessarily also flat, that is, metrics (\ref{m1})
form a flat pencil of metrics. The condition of
almost compatibility for metrics (\ref{m1})
has the form
\be
\eta^{sp} {\pa^2 \Phi \over \pa v^p \pa v^i}
{\pa^3 \Phi \over \pa v^s \pa v^j \pa v^k} =
\eta^{sp} {\pa^2 \Phi \over \pa v^p \pa v^k}
{\pa^3 \Phi \over \pa v^s \pa v^j \pa v^i} \label{m2}
\ee
and coincides with the condition of
compatible deformation of two Frobenius algebras
(this condition was derived and studied by the
present author in \cite{[36]}--\cite{mokh1}).
\et

In particular, in the present author's papers
\cite{[36]}--\cite{mokh1}, it is proved that
in the two-component case ($N=2$), for $\eta^{ij} =
\varepsilon^i \delta^{ij},$ $ \varepsilon^i = \pm 1,$
condition (\ref{m2}) is equivalent to the following
linear second-order partial differential equation with constant
coefficients:
\be
\alpha \left ( \varepsilon^1
{\pa^2 \Phi \over \pa (v^1)^2} -
 \varepsilon^2
{\pa^2 \Phi \over \pa (v^2)^2} \right ) = \beta
 {\pa^2 \Phi \over \pa v^1 \pa v^2},
\ee
where $\alpha$ and $\beta$ are arbitrary constants which are not
equal to zero simultaneously.

\section{Compatible metrics and the Nijenhuis tensor} \label{sect2}

Let us prove the second part of Theorem \ref{tmo1}.
In the previous section, it is proved, in particular,
that it always follows from compatibility (moreover,
even from almost compatibility)
of metrics that the corresponding Nijenhuis tensor vanishes
(Theorem \ref{pochti}).

Assume that a pair of metrics
$g^{ij}_1 (u)$ and $g^{ij}_2 (u)$ is nonsingular, that is,
the eigenvalues of this pair of metrics are distinct.
Furthermore, assume that the corresponding Nijenhuis tensor vanishes.
Let us prove that, in this case, the metrics
$g^{ij}_1 (u)$ and $g^{ij}_2 (u)$ are compatible
(their almost compatibility follows from Theorem \ref{pochti}).

It is obvious that the eigenvalues of the pair of metrics
$g^{ij}_1 (u)$ and $g^{ij}_2 (u)$  coincide with the
eigenvalues of the affinor $v^i_j (u)$. But it is well known
that if all eigenvalues of an affinor are distinct, then
it always follows from the vanishing of the Nijenhuis tensor
of this affinor that there exist local coordinates such that,
in these coordinates,
the affinor reduces to a diagonal form in the corresponding neighbourhood
\cite{nij1} (see also \cite{ha}).

So, further, we can consider that the affinor
$v^i_j (u)$ is diagonal in the local coordinates
 $u^1,...,u^N$, that is,
\be
 v^i_j (u) = \lambda^i (u) \delta^i_j,
\ee
where is no summation over the index $i$.
By assumption, the eigenvalues
$\lambda^i (u), \ i=1,...,N,$ coinciding with the
eigenvalues
of the pair of metrics $g^{ij}_1 (u)$ and $ g^{ij}_2 (u)$
are distinct:
\be
\lambda^i \neq \lambda^j  {\rm \ if\ } i\neq  j.
\ee

\bl              \label{n1}
If affinor $v^i_j (u)$ (\ref{aff}) is diagonal in certain local
coordinates and all its eigenvalues are distinct, then,
in these coordinates, the metrics $g^{ij}_1 (u) $ and
$g^{ij}_2 (u)$ are also necessarily diagonal.
\el

{\it Proof}.  Actually, we have
$$g^{ij}_1 (u) = \lambda^i (u) g^{ij}_2 (u).$$
It follows from symmetry of the metrics $g^{ij}_1 (u)$ and
$ g^{ij}_2 (u)$
that for any indices $i$ and  $j$
\be
(\lambda^i (u) - \lambda^j (u))  g^{ij}_2 (u) = 0,
\ee
where is no summation over indices, that is,
$$ g^{ij}_2 (u) = g^{ij}_1 (u) =0
{\rm \  if \ } i\neq  j.  {\ \rm \ q.e.d.}$$

\bl     \label{n2}
Let an affinor $w^i_j (u)$ be diagonal in certain
local coordinates
$u= (u^1,...,u^N)$, that is, $w^i_j (u) =\mu^i (u) \delta^i_j$.
\bitem
\item[1)] If all the eigenvalues $\mu^i (u), \ i=1,...,N,$
of the diagonal affinor are distinct,
that is,
$\mu^i (u) \neq \mu^j (u)$ for $i \neq  j$, then
the Nijenhuis tensor of this affinor vanishes if and only if
the $i$th eigenvalue $\mu^i (u)$ depends only on the coordinate $u^i.$
\item[2)] If all the eigenvalues coincide, then
the Nijenhuis tensor vanishes.
\item[3)] In the general case of an arbitrary diagonal affinor
$w^i_j (u) =\mu^i (u) \delta^i_j$,
the Nijenhuis tensor vanishes if and only if
\be
{\pa \mu^i \over \pa u^j} = 0
\ee
for all indices $i$ and $j$ such that $\mu^i (u) \neq \mu^j (u).$
\eitem
\el

{\it Proof}.
Actually, for any diagonal
affinor $w^i_j (u) =\mu^i (u) \delta^i_j,$
the Nijenhuis tensor $N^k_{ij} (u)$ has the form
$$N^k_{ij} (u) = (\mu^i - \mu^k) {\pa \mu^j \over \pa u^i} \delta^{kj}
- (\mu^j - \mu^k) {\pa \mu^i \over \pa u^j} \delta^{ki}$$
(no summation over indices).
Thus, the Nijenhuis tensor vanishes if and only if
for any indices $i$ and $j$
$$(\mu^i (u) - \mu^j (u)) {\pa \mu^i \over \pa u^j} =0,$$
where is no summation over indices. q.e.d.

It follows from Lemmas \ref{n1} and \ref{n2}
that for any nonsingular pair of
almost compatible metrics there always exist local coordinates
in which the metrics have the form
$$g^{ij}_2 (u) = g^i (u) \delta^{ij}, \ \ \
g^{ij}_1 (u) = \lambda^i (u^i) g^i (u) \delta^{ij},
\ \ \ \lambda^i = \lambda^i (u^i), \ i=1,...,N.$$

Moreover, we immediately derive that any pair of diagonal
metrics of the form $g^{ij}_2 (u) = g^i (u) \delta^{ij}$ and
$g^{ij}_1 (u) = f^i (u^i) g^i (u) \delta^{ij}$
for any nonzero functions $f^i (u^i),$ $ i=1,...,N,$
(here they can be, for example, coinciding nonzero constants, that is,
the pair of metrics may be ``singular'') is almost compatible,
since the corresponding Nijenhuis tensor always vanishes for any
pair of metrics of
this form.

We shall prove now that any pair of metrics of this
form is always compatible.
Then Theorems \ref{tmo1} and \ref{teomxx}
will be completely proved.

Consider two diagonal metrics of the form
$g^{ij}_2 (u) = g^i (u) \delta^{ij}$ and
$g^{ij}_1 (u) = f^i (u^i) g^i (u) \delta^{ij},$
where $f^i (u^i),$ $i=1,...,N,$ are arbitrary (possibly, complex)
nonzero functions of
single variables,
and consider their arbitrary linear combination
$$g^{ij} (u) = (\lambda_2 + \lambda_1 f^i (u^i)) g^i (u) \delta^{ij},$$
where $\lambda_1$ and $\lambda_2$ are arbitrary constants such that
$\det (g^{ij} (u)) \not\equiv 0.$

Let us prove that relation (\ref{kr}) is always fulfilled for
the corresponding tensors of Riemannian curvature.

Recall that for any diagonal metric
$\Gamma^i_{jk} (u) =0$ if all the indices $i, j, k$
are distinct.
Correspondingly,
$R^{ij}_{kl} (u) = 0$
if all the indices $i, j, k, l $  are distinct.
Besides, as a result of the well-known symmetries of
the tensor of Riemannian curvature, we have:
$$R^{ii}_{kl} (u) = R^{ij}_{kk} (u) =0,$$
$$R^{ij}_{il} (u) = -R^{ij}_{li} (u) = R^{ji}_{li} (u) =
- R^{ji}_{il} (u).$$

Thus, it is sufficient to prove relation (\ref{kr}) only
for the following components of the corresponding
tensors of Riemannian curvature:
 $R^{ij}_{il} (u)$,
where $i \neq j,$ $\ i \neq  l$.

For an arbitrary diagonal metric
$g^{ij}_2 (u) = g^i (u) \delta^{ij}$, we have
$$\Gamma^i_{2, ik} (u) = \Gamma^i_{2, ki} (u) =
- {1 \over 2 g^i (u)} {\pa g^i \over \pa u^k}, \ \ \
{\rm \ for\  any\ } i, k; $$
$$\Gamma^i_{2, jj} (u) = {1 \over 2} {g^i (u) \over
(g^j (u))^2 } {\pa g^j \over \pa u^i},\ \ \ i \neq j.$$

\bea
&&
R^{ij}_{2, il} (u) = g^i (u) R^j_{2, iil} (u) =
g^i (u) \left ( {\pa \Gamma^j_{2, il} \over \pa u^i}  -
{\pa \Gamma^j_{2, ii}  \over \pa u^l} +\right.
\\
&&
\left. \sum_{s=1}^N
\Gamma^j_{2, si} (u) \Gamma^s_{2, il} (u) -
\sum_{s=1}^N \Gamma^j_{2, sl} (u) \Gamma^s_{2, ii} (u) \right ).\nn
\eea

It is necessary to consider two the following
different cases separately:

1) $j \neq l$.

\bea
&&
R^{ij}_{2, il} (u) =
g^i (u) \left (
- {\pa \Gamma^j_{2, ii}  \over \pa u^l} +
\Gamma^j_{2, ii} (u) \Gamma^i_{2, il} (u) -\right. \label{r1xx}\\
&&
\left. \Gamma^j_{2, jl} (u) \Gamma^j_{2, ii} (u) -
\Gamma^j_{2, ll} (u) \Gamma^l_{2, ii} (u) \right )
=  \nn\\
&&
- {1 \over 2} g^i (u)  {\pa \over \pa u^l}
\left ( {g^j (u) \over (g^i (u))^2 } {\pa g^i \over \pa u^j }\right )
 - {1 \over 4} {g^j (u) \over (g^i (u))^2 }
{\pa g^i \over \pa u^j} {\pa g^i \over \pa u^l} +\nn\\
&&
{1 \over 4 g^i (u)}
{\pa g^i \over \pa u^j} {\pa g^j \over \pa u^l}
- {1 \over 4} {g^j (u) \over g^i (u) g^l (u) }
{\pa g^l \over \pa u^j} {\pa g^i \over \pa u^l}.\nn
\eea

Respectively, for the metric
$$g^{ij} (u) = (\lambda_2 + \lambda_1 f^i (u^i)) g^i (u) \delta^{ij},$$
we obtain (here we use that all the indices $i, j, l$ are distinct):
\bea
&&
R^{ij}_{il} (u)=
(\lambda_2 + \lambda_1 f^j (u^j)) \left [
- {1 \over 2} g^i (u)  {\pa \over \pa u^l}
\left ( {g^j (u) \over (g^i (u))^2 } {\pa g^i \over \pa u^j }\right )
 - \right. \label{rr1}\\
&&
\left. {1 \over 4} {g^j (u) \over (g^i (u))^2 }
{\pa g^i \over \pa u^j} {\pa g^i \over \pa u^l}
+ {1 \over 4 g^i (u)}
{\pa g^i \over \pa u^j} {\pa g^j \over \pa u^l}
- {1 \over 4} {g^j (u) \over g^i (u) g^l (u) }
{\pa g^l \over \pa u^j} {\pa g^i \over \pa u^l} \right ]= \nn\\
&&
\lambda_1 R^{ij}_{1, il} (u) + \lambda_2 R^{ij}_{2, il} (u). \nn
\eea

2) $j=l$.

\bea
&&
R^{ij}_{2, ij} (u) =
g^i (u) \left (
 {\pa \Gamma^j_{2, ij}  \over \pa u^i} -
{\pa \Gamma^j_{2, ii}  \over \pa u^j} +
\Gamma^j_{2, ii} (u) \Gamma^i_{2, ij} (u) + \right. \label{r2xx}\\
&&
\left. \Gamma^j_{2, ji} (u) \Gamma^j_{2, ij} (u)
-  \sum_{s=1}^N
\Gamma^j_{2, sj} (u) \Gamma^s_{2, ii} (u) \right )
=  \nn\\
&&
- {1 \over 2} g^i (u)  {\pa \over \pa u^i}
\left ( {1 \over g^j (u) } {\pa g^j \over \pa u^i }\right )
- {1 \over 2} g^i (u)  {\pa \over \pa u^j}
\left ( {g^j (u) \over (g^i (u))^2 } {\pa g^i \over \pa u^j }\right )
- \nn\\
&&
{1 \over 4} {g^j (u) \over (g^i (u))^2 }
{\pa g^i \over \pa u^j} {\pa g^i \over \pa u^j}
+ {1 \over 4} {g^i (u) \over (g^j(u))^2}
{\pa g^j \over \pa u^i} {\pa g^j \over \pa u^i}
- {1 \over 4 g^j (u)}
{\pa g^j \over \pa u^i} {\pa g^i \over \pa u^i} +\nn\\
&&
\sum_{s \neq  i} {1 \over 4} {g^s (u) \over g^i (u) g^j (u) }
{\pa g^j \over \pa u^s} {\pa g^i \over \pa u^s}.\nn
\eea

Respectively, for the metric
$$g^{ij} (u) = (\lambda_2 + \lambda_1 f^i (u^i)) g^i (u) \delta^{ij},$$
we obtain (here we use that the indices $i$ and $j$ are distinct):

\bea
&&
R^{ij}_{ij} (u) =
-{1 \over 2} (\lambda_2 + \lambda_1 f^i (u^i))
g^i (u)  {\pa \over \pa u^i}
\left ( {1 \over g^j (u) } {\pa g^j \over \pa u^i }\right )
-\\
&&
{1 \over 2} g^i (u)  {\pa \over \pa u^j}
\left ( {
(\lambda_2 + \lambda_1 f^j (u^j))
g^j (u) \over (g^i (u))^2 } {\pa g^i \over \pa u^j }\right )
-\nn\\
&&
{1 \over 4}
(\lambda_2 + \lambda_1 f^j (u^j))
{g^j (u) \over (g^i (u))^2 }
{\pa g^i \over \pa u^j} {\pa g^i \over \pa u^j}+ \nn\\
&&
{1 \over 4}
(\lambda_2 + \lambda_1 f^i (u^i))
 {g^i (u) \over (g^j(u))^2}
{\pa g^j \over \pa u^i} {\pa g^j \over \pa u^i}
-\nn\\
&&
{1 \over 4 g^j (u)}
{\pa g^j \over \pa u^i} {\pa
((\lambda_2 + \lambda_1 f^i (u^i))
g^i) \over \pa u^i} +
\nn\\
&&
{1 \over 4 g^i (u)}
{\pa g^i \over \pa u^j} {\pa
((\lambda_2 + \lambda_1 f^j (u^j))
g^j) \over \pa u^j} +
\nn\\
&&
\sum_{s \neq  i, \ s \neq  j}
{1 \over 4} {
(\lambda_2 + \lambda_1 f^s (u^s))
g^s (u) \over g^i (u) g^j (u) }
{\pa g^j \over \pa u^s} {\pa g^i \over \pa u^s}
=\nn\\
&&
\lambda_1 R^{ij}_{1, ij} (u) + \lambda_2 R^{ij}_{2, ij} (u).\nn
\eea

Theorems \ref{tmo1} and \ref{teomxx} are proved.
Thus, the complete explicit description
of nonsingular pairs of compatible and almost compatible
metrics is obtained.

\section{Equations for nonsingular pairs of \\
compatible flat metrics}

Now, let us consider, in detail, the problem on nonsingular pairs
of compatible flat metrics. It follows from
Theorem \ref{teomxx}
that it is
sufficient to classify all pairs of flat metrics of the
following special diagonal form
$g_2^{ij} (u) = g^i (u) \delta^{ij}$ and
$g_1^{ij} (u) = f^i (u^i) g^i (u) \delta^{ij},$
where $f^i (u^i),$ $i= 1,...,N,$ are arbitrary (possibly, complex)
functions of
single variables.

The  problem of description of diagonal flat metrics,
that is, flat metrics
$g_2^{ij} (u) = g^i (u) \delta^{ij},$
is a classical problem of differential geometry.
This problem is equivalent to the problem
of description of curvilinear orthogonal coordinate systems in
an N-dimensional
pseudo-Euclidean space and it was studied in detail and mainly solved
in the beginning of the 20th century (see \cite{da}).
Locally, such coordinate systems are determined by
$N(N-1)/2$ arbitrary functions of two variables.
Recently, Zakharov showed that the Lam\'{e}
equations describing curvilinear orthogonal coordinate systems
can be integrated by the inverse scattering method \cite{za}
(see also an algebraic-geometric approach in \cite{kri}).

The condition that the metric
$g_1^{ij} (u) = f^i (u^i) g^i (u) \delta^{ij}$ is also flat
exactly gives $N(N-1)/2$ additional equations linear with respect to
the functions $f^i (u^i)$.  Note that, in this case,
components (\ref{r1xx})
of the corresponding tensor of Riemannian curvature
automatically vanish as a result of formula (\ref{rr1}).
And the vanishing of components (\ref{r2xx}) gives
the corresponding $N(N-1)/2$ equations. In particular, in the case
$N=2$, this completely solves the problem of description
for nonsingular pairs of compatible two-component flat metrics.
In the next section, we present this complete description.
It is also very interesting to classify
all the $N$-orthogonal curvilinear coordinate systems in a
pseudo-Euclidean space (or, in other words, to classify the
corresponding functions $g^i (u)$) such that the functions
$f^i (u^i) = (u^i)^n$ define the corresponding
compatible flat metrics
(respectively, separately for $n=1$; $n=1, 2$; $n=1, 2, 3,$ and so on).

\bt \label{syst}
Any nonsingular pair of compatible flat metrics is described by
the following integrable nonlinear
system which is the special reduction of the Lam\'{e} equations:
\be
{\pa \beta_{ij} \over \pa u^k}
=\beta_{ik} \beta_{kj},\ \ \ i\neq j,\ \ i\neq k,\ \ j\neq k, \label{lam1}
\ee
\be
{\pa \beta_{ij} \over \pa u^i}+
{\pa \beta_{ji} \over \pa u^j}+\sum_{s\neq i,\
s\neq j} \beta_{si} \beta_{sj} =0,\ \ \ i\neq j, \label{lam2}
\ee
\bea
&&
\sqrt {f^i (u^i)} {\pa \left (\sqrt {f^i (u^i)} \beta_{ij}\right )
 \over \pa u^i}+
\sqrt {f^j (u^j)}
{\pa \left (\sqrt {f^j (u^j)}
\beta_{ji}\right ) \over \pa u^j}+ \label{lam3}\\
&&
\sum_{s\neq i,\
s\neq j} f^s (u^s) \beta_{si} \beta_{sj} =0,\ \ \ i\neq j, \nn
\eea
where $f^i (u^i),$ $i=1,...,N,$ are given arbitrary (possibly, complex)
functions
of single variables (these functions
are the eigenvalues of the pair of metrics).
\et

\br
Equations (\ref{lam1}) and (\ref{lam2}) are
the famous Lam\'{e} equations. Equations (\ref{lam3})
define a nontrivial nonlinear differential reduction of
the Lam\'{e} equations.
\er

{\it Proof}.
Consider the conditions of flatness for the diagonal metrics
$g^{ij}_2 (u) = g^i (u) \delta^{ij}$ and
$g^{ij}_1 (u) = f^i (u^i) g^i (u) \delta^{ij},$
where $f^i (u^i),$ $i=1,...,N,$ are arbitrary (possibly, complex)
functions of
the given single variables (but these functions
are not equal to zero identically).

As is shown in the previous section, for
any diagonal metric, it is sufficient to consider
the condition $R^{ij}_{kl} (u) =0$ (the condition of flatness for a
metric) only
for the following components of the tensor of Riemannian curvature:
 $R^{ij}_{il} (u)$,
where $i \neq j,$ $\ i \neq  l$.

Again as above,
for an arbitrary diagonal metric $g^{ij}_2 (u) = g^i (u) \delta^{ij}$,
it is necessary to consider two
the following different cases separately.

1) $j \neq l$.

\bea
&&
R^{ij}_{2, il} (u) = \label{r1}\\
&&
- {1 \over 2} g^i (u)  {\pa \over \pa u^l}
\left ( {g^j (u) \over (g^i (u))^2 } {\pa g^i \over \pa u^j }\right )
 - {1 \over 4} {g^j (u) \over (g^i (u))^2 }
{\pa g^i \over \pa u^j} {\pa g^i \over \pa u^l} + \nn\\
&&
{1 \over 4 g^i (u)}
{\pa g^i \over \pa u^j} {\pa g^j \over \pa u^l}
- {1 \over 4} {g^j (u) \over g^i (u) g^l (u) }
{\pa g^l \over \pa u^j} {\pa g^i \over \pa u^l}=0 .\nn
\eea

Introducing the standard classical notation
\bea
&&
g^i (u) = {1 \over (H_i (u))^2 },\ \ \ d\, s^2 = \sum_{i=1}^N
(H_i (u))^2 (d u^i)^2, \\
&&
\beta_{ik} (u) = {1 \over H_i (u)} {\pa H_k \over \pa u^i},\ \ \
i \neq k,
\eea
where $H_i (u)$ are the {\it Lam\'{e} coefficients} and
$\beta_{ik} (u)$ are the {\it rotation coefficients},
we derive that equations (\ref{r1}) are equivalent to
the equations
\be
{\pa^2 H_i \over \pa u^j \pa u^k} =
{1 \over H_j (u)} {\pa H_i \over \pa u^j} {\pa H_j \over \pa u^k}
+ {1 \over H_k (u)} {\pa H_k \over \pa u^j} {\pa H_i \over \pa u^k},
\label{h1}
\ee
where $i \neq j,$ $i \neq k,$ $j \neq k.$
Equations (\ref{h1}) are equivalent to equations (\ref{lam1}).

2) $j=l$.

\bea
&&
R^{ij}_{2, ij} (u) =\label{r2}\\
&&
- {1 \over 2} g^i (u)  {\pa \over \pa u^i}
\left ( {1 \over g^j (u) } {\pa g^j \over \pa u^i }\right )
- {1 \over 2} g^i (u)  {\pa \over \pa u^j}
\left ( {g^j (u) \over (g^i (u))^2 } {\pa g^i \over \pa u^j }\right )
- \nn\\
&&
{1 \over 4} {g^j (u) \over (g^i (u))^2 }
{\pa g^i \over \pa u^j} {\pa g^i \over \pa u^j}
+ {1 \over 4} {g^i (u) \over (g^j(u))^2}
{\pa g^j \over \pa u^i} {\pa g^j \over \pa u^i}
-\nn\\
&&
{1 \over 4 g^j (u)}
{\pa g^j \over \pa u^i} {\pa g^i \over \pa u^i} +
\sum_{s \neq  i} {1 \over 4} {g^s (u) \over g^i (u) g^j (u) }
{\pa g^j \over \pa u^s} {\pa g^i \over \pa u^s}=0.\nn
\eea

Equations (\ref{r2}) are equivalent to the equations
\bea
&&
{\pa \over \pa u^i} \left ( {1 \over H_i (u)}
{\pa H_j  \over \pa u^i} \right ) +
{\pa \over \pa u^j} \left ( {1 \over H_j (u)}
{\pa H_i  \over \pa u^j} \right ) + \label{h2}\\
&&
\sum_{s \neq i,\ s \neq j}
{1 \over (H_s (u))^2} {\pa H_i \over \pa u^s} {\pa H_j \over \pa u^s}=0,
\ \ \ i \neq j.\nn
\eea
Equations (\ref{h2}) are equivalent to equations (\ref{lam2}).

The condition that the metric
$g_1^{ij} (u) = f^i (u^i) g^i (u) \delta^{ij}$ is also flat
gives exactly $N(N-1)/2$ additional equations (\ref{lam3})
which are linear with respect to
the given functions $f^i (u^i)$.  Note that, in this case,
components (\ref{r1})
of the corresponding tensor of Riemannian curvature
automatically vanish.
And the vanishing of components (\ref{r2}) gives
the corresponding $N(N-1)/2$ additional equations.

Actually, for the metric
$g^{ij}_1 (u) = f^i (u^i) g^i (u) \delta^{ij},$
we have
\bea
&&
\widetilde H_i (u) = {H_i (u)  \over \sqrt {f^i (u^i)}},\\
&&
\widetilde \beta_{ik} (u) =
{1 \over \widetilde H_i (u)} {\pa \widetilde H_k \over
\pa u^i} =\\
&&
{\sqrt {f^i (u^i)} \over
\sqrt {f^k (u^k)}} \left ( {1 \over H_i (u)} {\pa H_k \over
\pa u^i}  \right ) =  {\sqrt {f^i (u^i)} \over
\sqrt {f^k (u^k)}}  \beta_{ik} (u), \ \ \ i \neq k. \nn
\eea
Respectively,
equations (\ref{lam1}) are also fulfilled for
the rotation coefficients $\widetilde \beta_{ik} (u)$
and equations (\ref{lam2}) for them give equations
(\ref{lam3}), which can be rewritten as follows
(as linear equations with respect to the functions $f^i(u^i)$):
\bea
&&
f^i (u^i) {\pa \beta_{ij} \over \pa u^i}+
{1 \over 2}  (f^i (u^i))' \beta_{ij} +
f^j (u^j) {\pa \beta_{ji} \over \pa u^j}+ \label{lam3x}\\
&&
{1 \over 2} (f^j (u^j))' \beta_{ji} +
\sum_{s\neq i,\
s\neq j} f^s (u^s) \beta_{si} \beta_{sj} =0,\ \ \ i\neq j. \nn
\eea
q.e.d.

\section{Two-component compatible flat metrics} \label{sect5}

Here we present the complete description for
nonsingular pairs of two-component compatible
flat metrics (see also \cite{[35]}, \cite{mokh1},
\cite{mokh2}, where an integrable
four-component nondiagonalizable homogeneous system of hydrodynamic type,
describing all the two-component compatible flat metrics, was derived and
investigated).

It is shown above that for any nonsingular pair
of two-component compatible metrics $g^{ij}_1 (u)$
and $g^{ij}_2 (u)$ there always exist local coordinates
$u^1,...,u^N$ such that
\be
(g^{ij}_2 (u) ) = \left ( \begin{array}{cc}
{\varepsilon^1 \over (b^1 (u))^2 } & 0\\
0 & {\varepsilon^2 \over (b^2 (u))^2 }
\end{array} \right ), \ \ \ \ \
(g^{ij}_1 (u) ) = \left ( \begin{array}{cc}
{\varepsilon^1  f^1 (u^1) \over (b^1 (u))^2 } & 0\\
0 & {\varepsilon^2  f^2 (u^2) \over (b^2 (u))^2 }
\end{array} \right ), \label{metr}
\ee
where $\varepsilon^i = \pm 1, \ i=1, 2;$
$b^i (u)$ and $f^i (u^i), \ i=1, 2,$ are arbitrary nonzero functions
of the corresponding variables.

\bl \label{moflat}
An arbitrary diagonal metric $g^{ij}_2 (u)$ (\ref{metr})
is flat if and only if the functions $b^i (u),$  $ i=1, 2,$
are solutions of the following linear system:
\be
{\pa b^2 \over \pa u^1} = \varepsilon^1
{\pa F \over \pa u^2} b^1 (u),\ \ \ \
{\pa b^1 \over \pa u^2} =
- \varepsilon^2 {\pa F \over \pa u^1} b^2 (u),
\label{sys}
\ee
where $F(u)$ is an arbitrary function.
\el

\bt
The metrics $g^{ij}_1 (u)$ and $g^{ij}_2 (u)$ (\ref{metr})
form a flat pencil of metrics if and only if
the functions $b^i (u), \ i=1, 2,$ are solutions of
the linear system (\ref{sys}), where the function $F(u)$ is
a solution of the following linear equation:
\be
\ \ \ 2 {\pa^2 F \over \pa u^1 \pa u^2} (f^1 (u^1) -  f^2 (u^2))
+ {\pa F \over \pa u^2} {d f^1 (u^1)  \over d u^1} -
{\pa F \over \pa u^1} {d f^2 (u^2) \over d u^2} =0. \label{lequa}
\ee
\et

If the eigenvalues of the pair of metrics
$g^{ij}_1 (u)$ and $g^{ij}_2 (u)$ are not only distinct but
also are not constants, then we can always choose local coordinates
such that $f^1 (u^1) = u^1,$  $f^2 (u^2) = u^2$ (see also the remark in
\cite{[34]}).
In this case, equation (\ref{lequa}) has the form
\be
2 {\pa^2 F \over \pa u^1 \pa u^2} (u^1 -  u^2)
+ {\pa F \over \pa u^2} -
{\pa F \over \pa u^1} = 0. \label{lequa2}
\ee

Let us continue this recurrent procedure for the metrics
$G^{ij}_{n+1} (u) = v^i_s (u) G^{sj}_n(u)$ with the help of
the affinor $v^i_j (u) = u^i \delta^i_j.$

\bt
Three metrics
\be
(G^{ij}_n (u) ) = \left ( \begin{array}{cc}
{\varepsilon^1  (u^1)^n \over (b^1 (u))^2 } & 0\\
0 & {\varepsilon^2   (u^2)^n \over (b^2 (u))^2 }
\end{array} \right ), \ \ \ \ \ n=0, 1, 2,\label{metr2}
\ee
form a flat pencil of metrics (pairwise compatible)
if and only if the functions
$b^i (u), \ i=1, 2,$ are solutions of linear system (\ref{sys}),
where
\be
F(u) = c \ln (u^1 - u^2),
\ee
$c$ is an arbitrary constant.
The metrics $G^{ij}_n (u), \ n=0,1,2,3,$ are flat only in the
most trivial case when $c =0$
and, respectively, $b^1 = b^1 (u^1),$ $b^2 = b^2 (u^2)$.

The metrics $G^{ij}_n (u), n=0,1,2,$ are flat and
the metric $G^{ij}_3 (u)$ is a metric of
nonzero constant Riemannian curvature $K \neq 0$
(in this case, the metrics $G^{ij}_n,$ $n=0, 1, 2, 3,$
form a pencil of metrics of constant Riemannian curvature)
if and only if
\be
(b^1 (u))^2 = (b^2 (u))^2 = {\varepsilon^2 \over 4 K}
(u^1 - u^2),\ \ \ \ \varepsilon^1 = - \varepsilon^2,
\ \ \ \ c = \pm {1 \over 2}.
\ee
\et

\section {Almost compatible metrics \\
that are not compatible} \label{counter}

\bl   \label{molem}
Two-component diagonal conformally Euclidean metric
$$g^{ij} (u) =  \exp (a (u)) \delta^{ij}, \
1 \leq i,j \leq 2,$$
is flat if and only if
the function $a(u)$ is harmonic, that is,
\be
\Delta a \equiv {\pa^2 a \over \pa (u^1)^2} +
{\pa^2 a \over \pa (u^2)^2} =0.
\ee
\el

In particular, the metric
$g_1^{ij} (u) =  \exp (u^1 u^2) \delta^{ij}, \
1 \leq i,j \leq 2$, is flat. It is obvious that the flat metrics
$g_1^{ij} (u) =  \exp (u^1 u^2) \delta^{ij}, \
1 \leq i,j \leq 2$, and $g_2^{ij} (u) = \delta^{ij}, \ 1 \leq i, j \leq 2,$
are almost compatible, the corresponding
Nijenhuis tensor (\ref{nij}) vanishes.
But it follows from Lemma \ref{molem} that these metrics are not
compatible, their sum is not a flat metric.

Similarly, it is also possible to construct other counterexamples to
Theorem \ref{tfer}. Moreover, the following statement is true.

\bp
Any nonconstant real harmonic function $a(u)$
defines a pair of almost compatible metrics
 $g_1^{ij} (u) =  \exp (a (u)) \delta^{ij}, \
1 \leq i,j \leq 2$, and $g_2^{ij} (u) = \delta^{ij}, \ 1 \leq i, j \leq 2,$
which are not compatible. These metrics are compatible if
and only if $a = a(u^1 \pm i u^2).$
\ep

Let us also construct almost compatible metrics of
constant Riemannian curvature that are not compatible.

\bl   \label{molem2}
Two-component diagonal conformally Euclidean metric
$$g^{ij} (u) =  \exp (a (u)) \delta^{ij}, \
1 \leq i,j \leq 2,$$
is a metric of constant Riemannian curvature $K$
if and only if the function $a(u)$ is a solution of the Liouville equation
\be
\Delta a \equiv {\pa^2 a \over \pa (u^1)^2} +
{\pa^2 a \over \pa (u^2)^2} = 2 K e^{- a (u)}.  \label{liuv}
\ee
\el

\bp
For the metrics
$g_1^{ij} (u) =  \exp (a (u)) \delta^{ij}, \
1 \leq i,j \leq 2$, and $g_2^{ij} (u) = \delta^{ij}, \ 1 \leq i, j \leq 2,$
the corresponding Nijenhuis tensor vanishes, that is, they are
always almost compatible.
But they are real compatible metrics of constant Riemannian curvature
$K$ and $0$, respectively,
only in the most trivial
case when the function $a (u)$ is constant and, consequently,
$K=0$.
Complex metrics are compatible if and only if
$a (u) = a (u^1 \pm i u^2)$ and, in this case, also $K=0$.
\ep

Note that all the one-component ``metrics'' are always
compatible, and all the one-component local Poisson structures
of hydrodynamic type are also always compatible.
Let us construct examples of
almost compatible metrics that are not compatible for any $N > 1$.

\bp
The metrics
$g_1^{ij} (u) =  b (u) \delta^{ij}, \
1 \leq i,j \leq N$ and
$g_2^{ij} (u) = \delta^{ij}, \ 1 \leq i, j \leq N,$
where $b (u)$ is an arbitrary function, are always almost compatible,
the corresponding Nijenhuis tensor vanishes.
But they are compatible real metrics only in the most trivial
case when the function $b(u)$ is constant.
Complex metrics are compatible if and only if either the function
 $b(u)$ is constant or  $N=2$ and $b(u) = b (u^1 \pm i u^2)$.
\ep

\section{Compatible flat metrics and\\
the Zakharov
method of \\ differential reductions}

Recall the Zakharov method for integrating
the Lam\'{e} equations (\ref{lam1}) and (\ref{lam2}) \cite{za}.

We must choose a matrix function $F_{ij} (s, s', u)$ and
solve the linear integral equation
\be
K_{ij}(s, s', u)  = F_{ij}(s, s', u) + \int_s^{\infty}
\sum_l K_{il} (s, q, u) F_{lj} (q, s', u) dq.  \label{int}
\ee
Then we obtain a one-parameter family of solutions of
the Lam\'{e} equations by the
formula
\be
\beta_{ij} (s, u) = K_{ji} (s, s, u).  \label{resh}
\ee

In particular,
if $F_{ij} (s, s', u) = f_{ij} (s-u^i, s'-u^j),$
where $f_{ij} (x, y)$ is an arbitrary matrix function of two variables,
then formula (\ref{resh}) produces solutions of equations (\ref{lam1}).
To satisfy equations (\ref{lam2}), Zakharov proposed to impose
on the ``dressing matrix function'' $F_{ij} (s-u^i, s' -u^j)$ a certain
additional linear differential relation. If $F_{ij} (s - u^i, s' - u^j)$
satisfy the Zakharov differential relation, then
the rotation coefficients $\beta_{ij} (u)$ satisfy additionally
equations
(\ref{lam2}).

Let us present a scheme for integrating all the system
(\ref{lam1})--(\ref{lam3}).

\bl
If both the function $F_{ij} (s-u^i, s' -u^j)$ and
the function
\be
\widetilde F_{ij} (s-u^i, s' -u^j) =
{\sqrt {f^j (u^j -s')} \over \sqrt {f^i (u^i -s)} }
F_{ij} (s-u^i, s' -u^j)   \label{f}
\ee
satisfy the Zakharov differential relation, then
the corresponding rotation coefficients $\beta_{ij} (u)$
(\ref{resh})
satisfy all the equations (\ref{lam1})--(\ref{lam3}).
\el

{\it Proof}.
Actually, if $K_{ij} (s, s', u)$ is the solution of
the linear integral equation (\ref{int}) corresponding to the
function $F_{ij} (s-u^i, s' -u^j)$, then
\be
\widetilde K_{ij} (s, s', u) =
{\sqrt {f^j (u^j -s')} \over \sqrt {f^i (u^i -s)} }
K_{ij} (s, s', u)
\ee
is the solution of (\ref{int}) corresponding to
function (\ref{f}).
It is easy to prove multiplying the integral equation
(\ref{int}) by
$$\sqrt {f^j (u^j -s')} \bigg /\sqrt {f^i (u^i -s)}.$$
The relation
\bea
&&
K_{ij}(s, s', u)  = F_{ij}(s -u^i, s' - u^j)
+ \label{int2x}\\
&&
\int_s^{\infty}
\sum_l K_{il} (s, q, u)
F_{lj} (q - u^l,  s'- u^j) dq  \nn
\eea
implies
\bea
&&
{\sqrt {f^j (u^j -s')} \over \sqrt {f^i (u^i -s)} }
K_{ij}(s, s', u)  =
{\sqrt {f^j (u^j -s')} \over \sqrt {f^i (u^i -s)} }
F_{ij}(s -u^i, s' - u^j)
+ \label{int2xx}\\
&&
\int_s^{\infty}
\sum_l
{\sqrt {f^l (u^l -q)} \over \sqrt {f^i (u^i -s)} }
K_{il} (s, q, u)
{\sqrt {f^j (u^j -s')} \over \sqrt {f^l (u^l -q)} }
F_{lj} (q - u^l,  s'- u^j) dq  \nn
\eea
and, finally, we have
\bea
&&
\widetilde K_{ij}(s, s', u)  = \widetilde F_{ij}(s -u^i, s' - u^j)
+ \label{int2}\\
&&
\int_s^{\infty}
\sum_l \widetilde K_{il} (s, q, u)
\widetilde F_{lj} (q - u^l,  s'- u^j) dq.  \nn
\eea

Then both $\widetilde \beta _{ij} (s, u) = \widetilde K_{ji} (s, s, u)$
and $\beta_{ij} (s, u) = K_{ji} (s, s, u)$
satisfy the Lam\'{e} equations (\ref{lam1}) and (\ref{lam2}).
Besides, we
have
\bea
&&
\widetilde \beta _{ij} (s, u) = \widetilde K_{ji} (s, s, u)
=\\
&&
{\sqrt {f^i (u^i -s)} \over \sqrt {f^j (u^j -s)} }
K_{ji} (s, s, u) =
{\sqrt {f^i (u^i -s)} \over \sqrt {f^j (u^j -s)} }
\beta_{ij} (s, u).\nn
\eea

Thus, in this case, the rotation coefficients $\beta_{ij} (u)$
exactly satisfy all the equations (\ref{lam1})--(\ref{lam3}), that is,
they generate the corresponding compatible flat metrics. q.e.d.

\section{Integrability of the equations for \\
nonsingular pairs of compatible flat metrics}

The Zakharov differential reduction can be written as follows \cite{za}:

\be
{\pa F_{ij} (s, s', u) \over \pa s'} +
{\pa F_{ji} (s', s, u) \over \pa s} = 0. \label{za1}
\ee

Thus,
to resolve these differential relations for the matrix
function $F_{ij} (s- u^i, s' -u^j)$,
we can introduce $N(N-1)/2$ arbitrary functions of two variables
$\Phi_{ij} (x, y),$ $i < j$, and put for $i<j$
\bea
&&
F_{ij} (s - u^i, s' - u^j) =
{\pa \Phi_{ij} (s-u^i, s' - u^j) \over \pa s},\label{f1}\\
&&
F_{ji} (s - u^i, s' - u^j) =
- {\pa \Phi_{ij} (s' -u^i, s - u^j) \over \pa s},\label{f1x}
\eea
and
\be
F_{ii} (s - u^i, s' - u^i) = {\pa \Phi_{ii} (s-u^i, s' - u^i)
\over \pa s}, \label{f2}
\ee
where $\Phi_{ii} (x, y)$, $i=1,...,N,$
are arbitrary skew-symmetric functions of two variables:
\be
\Phi_{ii} (x, y) = - \Phi_{ii} (y, x),
\ee
see \cite{za}.

For the function
\be
\widetilde F_{ij} (s-u^i, s' -u^j) =
{\sqrt {f^j (u^j -s')} \over \sqrt {f^i (u^i -s)} }
F_{ij} (s-u^i, s' -u^j),
\ee
the Zakharov differential relation (\ref{za1})
exactly gives $N(N-1)/2$ linear partial
differential equations of the second order
 for
$N(N-1)/2$  functions $\Phi_{ij} (s- u^i, s' - u^j), \ i < j,$
of two variables:
\bea
&&
{\pa \over \pa s'} \left (
{\sqrt {f^j (u^j -s')} \over \sqrt {f^i (u^i -s)} }
{\pa \Phi_{ij} (s-u^i, s' - u^j) \over \pa s}
\right ) -\\
&&
{\pa \over \pa s} \left (
{\sqrt {f^i (u^i -s)} \over \sqrt {f^j (u^j -s')} }
{\pa \Phi_{ij} (s-u^i, s' - u^j) \over \pa s'} \right )=0,
\ \ \ \ i < j,\nn
\eea
or, equivalently,
\bea
&&
2 {\pa^2 \Phi_{ij} (s - u^i, s' - u^j) \over
\pa u^i \pa u^j} \left (f^i (u^i -s) - f^j (u^j - s') \right )
+ \label{ss}\\
&&
{\partial \Phi_{ij} (s- u^i, s' -u^j) \over \pa u^j}
{d f^i (u^i-s) \over d u^i}  -\nn\\
&&
{\partial \Phi_{ij} (s- u^i, s' -u^j) \over \pa u^i}
{d f^j (u^j-s') \over d u^j} =0, \ \ \ \ i < j.\nn
\eea

It is very interesting that all these equations (\ref{ss})
for the functions $\Phi_{ij} (s- u^i, s' - u^j)$  are of the same
type as in the two-component case. In fact, these equations
coincide with the corresponding single equation (\ref{lequa})
for the two-component case.

Besides, for $N$  functions
$\Phi_{ii} (s-u^i, s' - u^i),$ we have also $N$ linear partial
differential
equations of the second order from the Zakharov differential relation
(\ref{za1}):
\bea
&&
{\pa \over \pa s'} \left (
{\sqrt {f^i (u^i -s')} \over \sqrt {f^i (u^i -s)} }
{\pa \Phi_{ii} (s-u^i, s' - u^i) \over \pa s}
\right ) +\\
&&
{\pa \over \pa s} \left (
{\sqrt {f^i (u^i -s)} \over \sqrt {f^i (u^i -s')} }
{\pa \Phi_{ii} (s' -u^i, s - u^i) \over \pa s'} \right ) = 0\nn
\eea
or, equivalently,
\bea
&&
2 {\pa^2 \Phi_{ii} (s - u^i, s' - u^i) \over
\pa s \pa s'} \left (f^i (u^i -s) - f^i (u^i - s') \right )
- \label{ss2}\\
&&
{\partial \Phi_{ii} (s- u^i, s' -u^i) \over \pa s}
{d f^i (u^i-s') \over d s'}  +\nn\\
&&
{\partial \Phi_{ii} (s- u^i, s' -u^i) \over \pa s'}
{d f^i (u^i -s) \over d s} =0.  \nn
\eea

Any solution of linear partial differential equations
(\ref{ss}) and (\ref{ss2}) generates a
one-parameter family of solutions of system
(\ref{lam1})--(\ref{lam3}) by linear relations and
formulas (\ref{f1}), (\ref{f1x}),
(\ref{f2}), (\ref{int}) and (\ref{resh}).
Thus, our problem is linearized.


\medskip

\begin{flushleft}
Centre for Nonlinear Studies,\\
L.D.Landau Institute for Theoretical Physics,\\
Russian Academy of Sciences,\\
ul. Kosygina, 2,\\
Moscow, 117940  Russia\\
e-mail: mokhov@genesis.mi.ras.ru, mokhov@landau.ac.ru\\

\smallskip

Department of Mathematics,\\
University of Paderborn,\\
Paderborn, Germany\\
e-mail: mokhov@uni-paderborn.de
\end{flushleft}

\end{document}